\documentclass[11pt]{article}%
\usepackage{amsmath}
\usepackage{amsfonts}
\usepackage{amssymb}
\usepackage{graphicx}%
\providecommand{\U}[1]{\protect\rule{.1in}{.1in}}
\providecommand{\U}[1]{\protect\rule{.1in}{.1in}}
\newcommand{\ignore}[1]{}
\newcounter{subscript}

\begin{document}

\title{An Application of Value Distribution Theory to Schr\"odinger Operators with
Absolutely Continuous Spectrum}
\date{}
\author{CHARLES FULTON\\
Department of Mathematical Sciences \\
Florida Insitute of Technology \\Melbourne, FL. 32901-6975 \\   \\
DAVID PEARSON \\
Department of Mathematics \\
University of Hull \\
Cottingham Road \\
\\Hull, HU6 7RX England \\United Kingsom \\  \\
STEVEN PRUESS \\
1133 N. Desert Deer Pass \\
Green Valley, Arizona 85614-5530 }
\maketitle


\abstract{We use the ``Value Distribution'' theory developed by Pearson and Breimesser to obtain a
sequence of functions in the eigenvalue parameter for some Sturm-Liouville problems which have the property of being ``uniformly asymptotically distributed''.}

\section{Introduction}

In this paper we apply a general result of Breimesser and Pearson
\cite{BP2003} to obtain a very general asymptotic result for Schr\"{o}dinger
equations with one regular endpoint and one Limit Point endpoint which
generates some absolutely continuous spectrum. The software package SLEDGE
\cite{SLEDGE,PRYCE}, developed by Pruess and Fulton, has as one of its options
the computation of the spectral function for Sturm-Liouville problems having
continuous spectrum. For cases where the spectrum is known to be absolutely
continuous, more recent research codes described in \cite{FPP2,FPP4,IMA} have
proven to be much better than the SLEDGE package. This paper will, in the
opinion of the authors, provide a backgound for development of even further
formulas and numerical algorithms for the spectral density function in cases
having some a.c. spectrum. The general theory of `value distribution' which is
providing the new approach described here has been developed by Pearson and
Breimesser \cite{DBP3,DBP4,BP2000,BP2003}.

For simplicity we restrict our attention to the problem for $x \in[a,\infty)$,
$a \in\mathrm{{R\!\!\!\!I} \;}$, with a Dirichlet boundary condition at $x=a$:%

\begin{align}
-y^{\prime\prime}+q(x)y  & =\lambda y,\qquad a\leq x<\infty,\label{1.1}\\
y(a)  & =0. \label{1.2}%
\end{align}

\noindent We pose the assumptions:

\begin{itemize}
\item \textbf{(i)} q(x) is continuous in $[a,\infty)$,

\item \textbf{(ii)} the \textbf{LP/O-N or LP/O} case occurs at $x=\infty$,

\item \textbf{(iii)} there exist one or more intervals where a.c. spectrum is
known to exist.
\end{itemize}

See Fulton and Pruess \cite[p. 114]{TOMS98} for these general endpoint
classifications: here \textbf{O} means equation (\ref{1.1}) is oscillatory at
$\infty$ for all real $\lambda$ and \textbf{O-N} means there exists
$\Lambda\in\mathrm{{R\!\!\!\!I} \;}$ such that equation (\ref{1.1}) is
oscillatory at $\infty$ for all $\lambda\in(\Lambda, \infty)$ and
nonoscillatory at $\infty$ for all $\lambda\in(-\infty, \Lambda)$. The
\textbf{LP/O-N} case includes potentials which satisfy
\[
\lim_{x\rightarrow\infty}q(x)=L
\]
for some finite limit ($\Lambda= L$ in this case) and periodic potentials
$q(x)$, as well as all quasi-periodic and almost periodic potentials. The
\textbf{LP/O} case includes potentials which are \textbf{LP} at $\infty$ and
satisfy
\[
\lim_{x\rightarrow\infty}q(x)=-\infty
\]
In all of the above mentioned cases, it is possible for a.c. spectrum to occur.

We now state a very general theorem on Herglotz functions of Breimesser and
Pearson \cite[Thm1,p.40]{BP2003} which will be applied to the Titchmarsh-Weyl
$m$-function associated with the above self-adjoint Sturm-Liouville problem
(\ref{1.1})-(\ref{1.2}) in the next section:

\textsc{Theorem 1}\hspace{0.05in}(Breimesser and Pearson, 2003) Let $f(z)$ be
an arbitrary Herglotz function. Let $S$ be an arbitrary Borel subset of
$\mathrm{{R\!\!\!\!I} \;}$, and let $\Lambda$ be any set of finite measure.
Define the angle $\theta$ subtended by $S$ at $z$ for all $z \in
\mathrm{{C\!\!\!I}}$ ($\mathrm{{C\!\!\!I}}$ = complex numbers) by
\begin{equation}
\theta(S,z)=\int_{S}\operatorname{Im}\left(  \frac{1}{t-z}\right)  dt.
\label{1.3}%
\end{equation}
Letting the boundary value of the Herglotz function on the real $\lambda$-axis
be
\begin{equation}
f^{+}(\lambda):=\lim_{\epsilon\rightarrow0}[f(\lambda+i\epsilon)] \label{1.4}%
\end{equation}
we define,
\begin{equation}
\omega(\lambda,S,f^{+}):=\frac{1}{\pi}\theta(S,f^{+}), \label{1.5}%
\end{equation}
and
\begin{equation}
\omega(\lambda,S,f(\lambda+i\epsilon)):=\frac{1}{\pi}\theta(S,f(\lambda
+i\epsilon)). \label{1.6}%
\end{equation}
Then we have
\begin{align}
&  \left\vert \int_{\Lambda}\omega(\lambda,S,f^{+})\,d\lambda-\int_{\Lambda
}\omega(\lambda,S,f(\lambda+i\epsilon)\,d\lambda\right\vert \label{1.7}\\
&  \leq\frac{1}{\pi}\int_{\Lambda}\theta(\Lambda^{c},\lambda+i\epsilon
)\,d\lambda\longrightarrow0,\,as\;\epsilon\rightarrow0.\nonumber
\end{align}
where $\theta(\Lambda^{C}, \lambda+ i \epsilon)$ is the angle subtended at
$\lambda+i\epsilon$ by $\Lambda^{C} = \mathrm{{R\!\!\!\!I} \; \backslash
\Lambda}$.

\section{A Uniformly Asymptotically Distributed Sequence of Functions}

\setcounter{equation}{0}

We define a fundamental system $\{ u, v \}$ of solutions of the Schr\"odinger
equation (\ref{1.1}) at $x=a$ by the initial conditions
\begin{equation}
\left[
\begin{array}
[c]{ll}%
u(a,\lambda) & v(a,\lambda)\\
u^{\prime}(a,\lambda) & v^{\prime}(a,\lambda)
\end{array}
\right]  =\left[
\begin{array}
[c]{cc}%
1 & 0\\
0 & 1
\end{array}
\right]  , \lambda\in\mathrm{{C\!\!\!I}. }\label{2.1}%
\end{equation}
The Titchmarsh-Weyl $m-$function for the problem (\ref{1.1})-(\ref{1.2}) is
defined (for Im $z \neq0$) by
\begin{equation}
u(x,z)+m(z)v(x,z) \in L_{2}(a,\infty). \label{2.2}%
\end{equation}
The boundary value of the $m$-function on the real $\lambda$-axis is then
defined by
\begin{equation}
m^{+}(\lambda):=\lim_{\epsilon\rightarrow0}[m(\lambda+i\epsilon)]=A(\lambda
)+iB(\lambda), \label{2.3}%
\end{equation}
where we have introduced the real-valued functions $A(\lambda)$, $B(\lambda)$
as the real and imaginary parts of $m^{+}(\lambda)$.

The standard method of approximating spectral quantities for the half line
problem (\ref{1.1})-(\ref{1.2}) is to make use of spectral quantities
associated with the truncated regular problem on $[a,b]$ for $b \in(a,
\infty)$:%

\begin{align}
-y^{\prime\prime}+q(x)y  & =\lambda y,\qquad a\leq x\leq b\label{2.4}\\
y(a)  & =0\label{2.5}\\
y(b)  & =0\label{2.6}%
\end{align}
We let $m_{b}(z)$ be the $m$-function associated with the Sturm-Liouville
problem (\ref{2.4})-(\ref{2.6}), that is, for all $z$ not eigenvalues,
\[
u(b,z)+m_{b}(z)v(b,z)=0,
\]
\noindent or
\begin{equation}
m_{b}(z):=-\frac{u(b,z)}{v(b,z)}. \label{2.7}%
\end{equation}
The boundary value of this meromorphic $m_{b}$ function on the real $\lambda
$-axis is then well-defined for all $\lambda$ not eigenvalues of
(\ref{2.4})-(\ref{2.6}) as
\begin{equation}
m_{b}^{+}:=\lim_{\epsilon\rightarrow0}[m_{b}(\lambda+i\epsilon)]=-\frac
{u(b,\lambda)}{v(b,\lambda)}. \label{2.8}%
\end{equation}
Here, $m_{b}^{+}(\lambda)$ is necessarily real-valued.

As in Theorem 1, we let $S$ be any Borel subset of $\mathrm{{R\!\!\!\!I} \;}$
and $\Lambda$ be any finite interval in $\mathrm{{R\!\!\!\!I} \;}$. For $x
\in[a, \infty)$ and all $\lambda\in\mathrm{{R\!\!\!\!I} \;}$, not eigenvalues
(with $b=x$) of (\ref{2.4})-(\ref{2.6}), we put
\begin{equation}
F(x,\lambda):=-\frac{u(x,\lambda)}{v(x,\lambda)}=m_{x}^{+}(\lambda),
\label{2.9}%
\end{equation}
and define
\begin{equation}
F_{x}^{-1}(S)=\{\lambda\in\mathrm{{R\!\!\!\!I} \; \,| \;F(x,\lambda)\in S\}.
}\label{2.10}%
\end{equation}
For any real number $X$ it follows from (\ref{1.3}) that the angle subtended
by $S$ at $X$ is
\begin{equation}
\theta(X,S)=\left\{
\begin{array}
[c]{c}%
\pi\;\;if\;X\in S\\
0\;\;if\;X\notin S.
\end{array}
\right.  \label{2.11}%
\end{equation}
Hence for $X = F(x, \lambda) = m_{x}^{+}(\lambda) \in\mathrm{{R\!\!\!\!I} \;}$
we have
\begin{equation}
\omega(\lambda,S,m_{x}^{+}):=\frac{1}{\pi}\theta(S,F(x,\lambda))=\chi
(F_{x}^{-1}(S)) \label{2.12}%
\end{equation}
where $\chi(F_{x}^{-1}(S))$ is the characteristic function of the set
$F_{x}^{-1}(S)$ on the real $\lambda$-axis. Hence,
\begin{equation}
\int_{\Lambda}\omega(\lambda,S,m_{x}^{+})\,d\lambda=\int_{\Lambda}\chi
(F_{x}^{-1}(S))\,d\lambda=\mu\left(  \Lambda\cap F_{x}^{-1}(S)\right) ,
\label{2.13}%
\end{equation}
where $\mu(\Lambda\cup F_{x}^{-1}(S))$ denotes the Lebesgue measure of
$\Lambda\cap F_{x}^{-1}(S)$. We then have the following theorem:

\textsc{Theorem 2}. (i) For any Borel set $S$ and any finite interval
$\Lambda$ we have
\begin{equation}
\lim_{x\rightarrow\infty}\mu\left(  \Lambda\cap F_{x}^{-1}(S)\right)
=\int_{\Lambda}\omega(\lambda,S,m^{+})\,d\lambda\label{2.14}%
\end{equation}
where $\omega(\lambda, S, m^{+})$ is defined as in (\ref{1.5}) for the
boundary value (\ref{2.3}) of the Titchmarsh-Weyl $m$-function for the half
line problem (\ref{1.1})-(\ref{1.2}).

\noindent(ii) Let $S_{\lambda} = (\alpha(\lambda), \beta(\lambda))$, where
$\alpha(\lambda)$ and $\beta(\lambda)$ are bounded, Lebesgue measurable
functions on the interval $\Lambda$. Then we have%

\begin{equation}
\lim_{x\rightarrow\infty}\mu\left(  \Lambda\cap F_{x}^{-1}(S_{\lambda
})\right)  =\int_{\Lambda}\omega(\lambda,S_{\lambda},m^{+})\,d\lambda.
\label{2.15}%
\end{equation}


\vskip 4pt \noindent Proof of (i): Replacing the left-hand side of
(\ref{2.14}) by the left-hand side of (\ref{2.13}) we observe that
\begin{align*}
&  \left\vert \int_{\Lambda}\,\mu\left(  \Lambda\cap F_{x}^{-1}(S)\right)
\,d\lambda-\int_{\Lambda}\omega(\lambda,S,m^{+})\,d\lambda\right\vert \\
&  =\left\vert \int_{\Lambda}\omega(\lambda,S,m_{x}^{+})\,d\lambda
-\int_{\Lambda}\omega(\lambda,S,m^{+})\,d\lambda\right\vert \\
&  \leq\left\vert \int_{\Lambda}(\omega(\lambda,S,m_{x}^{+})\,-\omega
(\lambda,S,m_{x}(\lambda+i\epsilon))\,d\lambda\right\vert \\
&  +\left\vert \int_{\Lambda}(\omega(\lambda,S,m_{x}(\lambda+i\epsilon
)\,-\omega(\lambda,S,m(\lambda+i\epsilon))\,d\lambda\right\vert \\
&  +\left\vert \int_{\Lambda}(\omega(\lambda,S,m_{x}(\lambda+i\epsilon
)\,-\omega(\lambda,S,m^{+}))\,d\lambda\right\vert .
\end{align*}
Since $m_{x}(z)$ and $m(z)$ are Herglotz functions (by virtue of being
Titchmarsh-Weyl $m$-functions for the problems (\ref{1.1})-(\ref{1.2}) and
(\ref{2.4})-(\ref{2.6}), Theorem 1 applies to give the limit of the first and
third terms to be zero as $\epsilon\downarrow0$. For the second term we have
for the integrand (for fixed $\epsilon> 0$), as $x \to\infty$, that:
\begin{align*}
&  \left\vert \omega(\lambda,S,m_{x}(\lambda+i\epsilon))\,-\omega
(\lambda,S,m(\lambda+i\epsilon))\right\vert \\
&  =\left\vert \int_{S}\,\left[  \operatorname{Im}\left(  \frac{1}%
{t-m_{x}(\lambda+i\epsilon)}\right)  -\operatorname{Im}\left(  \frac
{1}{t-m(\lambda+i\epsilon)}\right)  \right]  dt\right\vert \\
&  =\int_{S}\,\operatorname{Im}\left(  \frac{-m(\lambda+i\epsilon
)+m_{x}(\lambda+i\epsilon)}{(t-m_{x}(\lambda+i\epsilon))(t-m(\lambda
+i\epsilon))}\right)  dt\\
&  \leq r_{x}(\lambda+i\epsilon)\cdot\int_{S}\,\left\vert \frac{1}%
{(t-m_{x}(\lambda+i\epsilon))(t-m(\lambda+i\epsilon))}\right\vert \,dt\\
&  \leq r_{x}(\lambda+i\epsilon)\cdot\left(  \int_{S}\,\left\vert \frac
{1}{(t-m(\lambda+i\epsilon))}\right\vert ^{2}\,dt\,+1\right) \\
&  \leq K_{\epsilon}\,\left(  \max_{\lambda\in\Lambda}\,r_{x}(\lambda
+i\epsilon)\right)  ,
\end{align*}
where $K_{\epsilon}$ is independent of $\lambda\in\Lambda$, and $r_{x}%
(\lambda+ i \epsilon)$ is the radius of the circle $C_{x, \lambda+i\epsilon}$
in the $m$-plane which contracts to the limit point $m(\lambda+ i \epsilon)$.
Hence, for fixed $\epsilon>0$,
\[
\lim_{x\rightarrow\infty}\omega(\lambda,S,m_{x}(\lambda+i\epsilon
)\,)=\omega(\lambda,S,m(\lambda+i\epsilon))
\]
uniformly for $\lambda\in\Lambda$. It follows that the above three terms can
be made arbitrarily small by fixing $\epsilon> 0$ sufficiently small and then
choosing $x$ sufficiently large. \quad\rule{2.2mm}{3.2mm}

\vskip 4pt \noindent Proof of (ii). First we assume that $S_{\lambda}=
(\alpha(\lambda), \beta(\lambda))$, $\alpha(\lambda) < \beta(\lambda)$ for all
$\lambda\in\Lambda$, where $\alpha(\lambda)$ and $\beta(\lambda)$ are step
functions on $\Lambda$, each taking only finitely many values. Letting
$\{\alpha_{j} \}, j = 1, \ldots, N$, and $\{\beta_{k} \}, k=1, \ldots, M$,
denote the values which $\alpha(\lambda)$ and $\beta(\lambda)$ may have for
$\lambda\in\Lambda$, we define
\[
C_{j,k} = \{ \lambda\in\Lambda\: | \:\alpha(\lambda) =\alpha_{j} \mbox{ and }
\beta(\lambda) = \beta_{k} \}.
\]
Then the $C_{j,k}$ are disjoint sets and $\lambda\in C_{j,k} \implies
\alpha_{j} < \beta_{k}$. Putting $S_{j,k} := (\alpha_{j}, \beta_{k})$ whenever
$\alpha_{j} < \beta_{k}$, we define for $x \in(a, \infty)$,
\[
F_{x}^{-1}(S_{j,k})=\left\{  \lambda\in\mathrm{{R\!\!\!\!I} \; \, | \;
F(x,\lambda)\in S_{j,k}}\right\} .
\]
Then the statement (\ref{2.14}) from part(i) applies with $S = S_{j,k}$, that
is,
\[
\lim_{x\rightarrow\infty}\mu\left(  \Lambda\cap F_{x}^{-1}(S_{j,k})\right)
=\int_{\Lambda}\omega(\lambda,S_{j,k},m^{+})\,d\lambda.
\]
Since the $C_{j,k}$ are disjoint sets on the $\lambda$-axis, it follows that
\begin{align*}
\mu\left(  \Lambda\cap F_{x}^{-1}(S_{\lambda})\right)   &  =\sum_{j,k}%
\mu\left(  \Lambda\cap C_{j,k}\cap F_{x}^{-1}(S_{j,k})\right) \\
&  \longrightarrow\sum_{j,k}\int_{\Lambda}\omega(\lambda,S_{j,k}%
,m^{+})\,d\lambda\\
&  =\sum_{j,k}\int_{\Lambda}\int_{\alpha_{j}}^{\beta_{k}}\,\left[
\operatorname{Im}\left(  \frac{1}{t-m^{+}(\lambda)}\right)  \right]
dt\,\;d\lambda\\
&  =\int_{\Lambda}\int_{\alpha(\lambda)}^{\beta(\lambda)}\, \left[
\operatorname{Im}\left(  \frac{1}{t-m^{+}(\lambda)}\right)  \right]
dt\,\;d\lambda\\
&  =\int_{\Lambda}\omega(\lambda,S_{\lambda},m^{+})\,d\lambda
,\;as\;x\longrightarrow\infty,
\end{align*}
which proves the statement (\ref{2.15}) when $S_{\lambda} = (\alpha(\lambda),
\beta(\lambda))$ and $\alpha(\lambda)$ and $\beta(\lambda)$ are simple
functions on $\Lambda$. The case where $\alpha(\lambda)$ and $\beta(\lambda)$
are general bounded Lebesgue measurable functions can then be handled by
approximating $\alpha(\lambda)$ and $\beta(\lambda)$ above and below by simple
functions. \quad\rule{2.2mm}{3.2mm}

Similarly to the analysis in Fulton, Pearson and Pruess \cite{IMA} and
Al-Naggar and Pearson \cite{AlN2}, we now introduce the space of quadratic
forms of solutions of equation (\ref{1.1}): For $\lambda\in
\mathrm{{R\!\!\!\!I} \;}$, let
\begin{equation}
\mathcal{R} = \left\{  R(x, \lambda) = au^{2} + b u v + c v^{2} \, | \, a, b,
c \in\mathrm{{R\!\!\!\!I} }\; \right\} \label{2.16}%
\end{equation}
be the vector space of all quadratic forms of solutions of equation
(\ref{1.1}). If $y_{1} = c_{1}u + d_{1} v$ and $y_{2} = c_{2}u + d_{2} v$, it
is readily verified that $y_{1}^{2}$, $y_{1}y_{2}$, $y_{2}^{2}$ also belong to
$\mathcal{R}$. Hence $\mathcal{R}$ is a 3-dimensional space with $\{u^{2}, uv,
v^{2}\}$ as a basis. Also associated with $\mathcal{R}$ is the solution space
of the first order Appell system (see \cite{IMA,APPELL}) of equations
\begin{equation}
\frac{\displaystyle dU}{\displaystyle dx}=\frac{d}{dx}\left[
\begin{array}
[c]{c}%
P\\
Q\\
R
\end{array}
\right]  =\left[
\begin{array}
[c]{ccc}%
0 & \lambda-q & 0\\
-2 & 0 & 2(\lambda-q)\\
0 & -1 & 0
\end{array}
\right]  \cdot\left[  \overset{}{%
\begin{array}
[c]{c}%
P\\
Q\\
R
\end{array}
}\right] , \label{2.17}%
\end{equation}
which was used extensively in \cite{FPP2,FPP4,IMA} and turns out to be a
useful companion system to the Schr\"odinger equation (\ref{1.1}) in the study
of cases involving some a.c. spectrum. One choice of basis for the solution
space of (\ref{2.17}) is
\begin{equation}
\left[  U_{1},U_{2},U_{3}\right]  =\left[
\begin{array}
[c]{ccc}%
(u^{\prime})^{2} & u^{\prime}v^{\prime} & (v^{\prime})^{2}\\
-2uu^{\prime} & -[u^{\prime}v+uv^{\prime}] & -2vv^{\prime}\\
u^{2} & uv & v^{2}%
\end{array}
\right] . \label{2.18}%
\end{equation}
Also on the solution space of (\ref{2.17}) we can define an indefinite inner
product by
\begin{equation}
\langle U,\tilde U\rangle:=2(P \tilde R + \tilde P R )-Q \tilde Q=const,
\mbox{ independent of } x\in\lbrack a,\infty), \label{2.19}%
\end{equation}
and the above solutions $\{U_{1}, U_{2}, U_{3} \}$ are orthogonal with respect
to this inner product. If
\[
U(x,\lambda)=a(\lambda)U_{1}(x,\lambda)+b(\lambda)U_{2}(x,\lambda
)+c(\lambda)U_{3}(x,\lambda),
\]
\[
\tilde U(x,\lambda)=\tilde a(\lambda)U_{1}(x,\lambda)+\tilde b(\lambda
)U_{2}(x,\lambda)+  \tilde c(\lambda)U_{3}(x,\lambda),
\]
then this inner product can also be expressed as
\begin{equation}
\langle U, \tilde U \rangle= 2 (a \tilde c + c \tilde a) - b \tilde
b.\label{2.20}%
\end{equation}

A solution $U$, or equivalently, its third component $R$ is said to be
normalized with respect to this indefinite inner product (compare \cite[Equa
4.10]{IMA} and \cite[Equa 7]{AlN2}) if
\begin{equation}
\langle U,U\rangle:=4PR-Q^{2}=2RR^{\prime\prime\prime})^{2}-4(q-\lambda
)R^{2}=4ac-b^{2}=4. \label{2.21}%
\end{equation}
We call $4ac - b^{2}$ the ``discriminant" of the quadratic form $\mathcal{R} =
au^{2} + buv + c v^{2}$. Also, it is readily shown from (\ref{2.17}) that the
third component of $R$ of any solution $U$ of (\ref{2.17}) must satisfy the
third order linear equation
\begin{equation}
R^{^{\prime\prime\prime}}+4(\lambda-q)R^{\prime}-2q^{\prime}R=0.\label{2.22}%
\end{equation}
The solution space of (\ref{2.22}) is the quadratic form space $\mathcal{R}$
and so there is a one-to-one correspondence between the solution space of
(\ref{2.22}) and (\ref{2.17}). Indeed, for every element, $R(x, \lambda)$, of
the quadratic form space $\mathcal{R}$ we can generate from the Appell
equations (\ref{2.17}) unique functions $P(x, \lambda)$ and $Q(x, \lambda)$
such that $(P,Q,R)^{T}$ is a solution of (\ref{2.17}).

We now cite the following lemma (see Al-Naggar and Pearson \cite[Lemma 1 and
Theorem 2]{AlN1}) which gives a sufficient condition for existence of a.c.
spectrum in a given interval $I \subset(-\infty, \infty)$:

\vskip6pt \noindent\textsc{Definition.} The Sturm-Liouville equation
(\ref{1.1}) satisfies \textit{Condition A} for a given real value of $\lambda$
if and only if there exists a complex-valued solution $y(x,\lambda)$ of
(\ref{1.1}) \ for which
\begin{equation}
\lim_{N\rightarrow\infty}\frac{\displaystyle\int_{0}^{N}y(x,\lambda)^{2}%
\,dx}{\displaystyle\int_{0}^{N}|y(x,\lambda)|^{2}\,dx}=0. \label{2.23}%
\end{equation}

\vskip 6pt \noindent\textsc{Lemma } (\textit{Al-Naggar and Pearson}). Let $I
\subset(-\infty,\infty)$ be an interval on which Condition A holds for the
general equation (\ref{1.1}), and let the fundamental system $\{u(\cdot
,\lambda),v(\cdot,\lambda)\}$ be defined by the initial conditions (\ref{2.1})
at $x = a$. Then \newline\textbf{(i)} There exists a complex valued function
$M(\lambda)$ on $I$ which is uniquely defined by the properties:
\begin{align}
&  (a) \text{ \ \ \ \ } \operatorname{Im}[M(\lambda)] > 0 \text{ \ \ and
\ \ }(b) \text{ \ \ } \lim_{N \rightarrow\infty} \frac{\displaystyle \int%
_{x_{0}}^{N} (u(x,\lambda)+ M(\lambda) v(x,\lambda))^{2} \, dx.}%
{\displaystyle \int_{x_{0}}^{N} |(u(x,\lambda)+ M(\lambda) v(x,\lambda)|^{2}
\, dx} = 0.\label{2.24}%
\end{align}
\textbf{(ii)} For $\lambda\in I$ the function $M$ in \textbf{(i)} is the
boundary value of the Titchmarsh-Weyl $m$-fuction defined by (\ref{2.3}), that
is,
\begin{equation}
\label{2.25}\mathit{M(\lambda)= \lim_{\epsilon\downarrow0}} [ m(\lambda+
i\epsilon)] = A(\lambda)+iB(\lambda).
\end{equation}
\textbf{(iii)} $I \subset\sigma_{ac}$, where $\sigma_{ac}$ is the absolutely
continuous spectrum of the Sturm-Liouville problem (\ref{1.1})-(\ref{1.2}),
and on $I$ the spectral density function is given by
\begin{equation}
f(\lambda):=\frac{1}{\pi}\lim_{\epsilon\rightarrow0}\,\left(
\operatorname{Im}[m(\lambda+i\epsilon)]\right)  =\frac{\operatorname{Im}%
[m^{+}(\lambda)]}{\pi}=\frac{B(\lambda)}{\pi}. \label{2.26}%
\end{equation}

In accordance with the above lemma, we can formulate the assumption (iii) of
a.c. spectrum as follows: \vskip 4pt \noindent(iii) Let $I \subset(-\infty,
\infty)$ be an interval on which $B(\lambda) = \operatorname{Im}[m^{+}%
(\lambda)] > 0$. \vskip 4pt \noindent Under this assumption, we now introduce
a polar coordinate representation of $u+m^{+}v = (u+Av) + i Bv$, namely
\begin{align}
u+Av  &  =r_{0}(x,\lambda)\cos(\theta_{0}(x,\lambda))\label{2.27}\\
Bv  &  =r_{0}(x,\lambda)\sin(\theta_{0}(x,\lambda))\nonumber
\end{align}
so that $\theta_{0}(x,\lambda)$ may be defined for $\lambda\in I$ and $x
\in[a, \infty)$ by
\begin{equation}
\cot(\theta_{0}(x,\lambda))=\left(  \frac{u(x,\lambda)+A(\lambda)v(x,\lambda
)}{B(\lambda)v(x,\lambda)}\right) . \label{2.28}%
\end{equation}
Taking the inverse of the cotangent, we have
\begin{equation}
\theta_{0}(x,\lambda)=\cot^{-1}\left(  \frac{u(x,\lambda)+A(\lambda
)v(x,\lambda)}{B(\lambda)v(x,\lambda)}\right)  \label{2.29}%
\end{equation}
from which we deduce that
\begin{equation}
\theta_{0}^{\prime}(x,\lambda)=\frac{B}{\left(  u+Av\right)  ^{2}+B^{2}v^{2}%
}>0 \label{2.30}%
\end{equation}
for all $\lambda\in I$ and all $x \in[a,\infty)$. This suggests that we look
at the quadratic form defined by
\begin{align}
R_{0}(x, \lambda)  &  := \frac{|u + m^{+}(\lambda) v|^{2}}{\mbox{Im }
m^{+}(\lambda)}\label{2.31}\\
&  = \frac{|u + (A + iB)v|^{2}}{B}\nonumber\\
&  = \frac{(u + Av)^{2} + B^{2} v^{2}}{B}\nonumber\\
&  = \frac{1}{B} u^{2} + \frac{2A}{B} uv + \frac{A^{2} +B^{2}}{B}
v^{2}.\nonumber
\end{align}
From (\ref{2.30}) and (\ref{2.31}) we see that the monotone increasing
function $\theta_{0}(x, \lambda)$, $x \to\infty$, is related to $R_{0}(x,
\lambda)$ by
\begin{equation}
\theta_{0}(x,\lambda)=\int_{a}^{x}\;\frac{1}{R_{0}(t,\lambda)}dt=\int_{a}%
^{x}\,\frac{B}{\left(  u+Av\right)  ^{2}+B^{2}v^{2}}\,dt. \label{2.32}%
\end{equation}
To keep $\theta_{0}(x, \lambda)$ lying in $[0, \pi]$ for all $x \in[a,
\infty)$, we restrict it by setting it back to zero whenever $x$ passes
through a zero of $v(x, \lambda)$, that is, we set
\begin{equation}
\tilde\theta_{0}(x,\lambda):=\theta_{0}(x,\lambda)|_{\operatorname{mod}\pi}
:=\theta_{0}(x,\lambda)-n\pi\;for\;x\in\;(x_{n},x_{n+1}), \label{2.33}%
\end{equation}
where $x_{0} = a$, $x_{n} = x_{n}(\lambda) = n\mbox{th zero of }v$. We can now
prove the following theorem concerning the asymptotic behavior of the sequence
of functions of $\lambda$, $\{\tilde\theta_{0}(x,\lambda) \} $ as $x \to
\infty$.

\vskip 8pt \noindent\textsc{Theorem 3}. Let  I  be an interval on which  $B(\lambda) > 0$.
Let $\Lambda \subset I$ be a bounded interval
and assume $C(\lambda)$ and $D(\lambda)$ are measurable functions on $\Lambda$
satisfying for all $\lambda\in\Lambda$
\[
0 \leq C(\lambda) < D(\lambda) \leq\pi.
\]
Then the set of functions of $\lambda$, $\{\tilde\theta_{0}(x,\cdot):I\rightarrow[0,\pi] | x \in (a,\infty)\}$
has the asymptotic
property as $x \to\infty$
\begin{equation}
\lim_{x \to\infty} \mu( \Lambda\cap\{\lambda\:| \; C(\lambda) < \tilde
\theta_{0}(x, \lambda) < D(\lambda) \} = \frac{1}{\pi} \int_{\Lambda} \; \{
D(\lambda) - C(\lambda)\} d \lambda.\label{2.34}%
\end{equation}
\vskip 4pt \noindent\textsc{Proof}. For $y = \cot\beta$ with $\beta\in(0,
\pi)$, we define $\beta(y) = cot^{-1}(y)$ mapping $(-\infty, \infty)$ onto
$(0, \pi)$. Since $\cot\beta$ is monotonically descreasing from $\beta= 0$ to
$\beta= \pi$, it follows that for $x \in(x_{n}, x_{n+1})$,
\[
\cot(C(\lambda))>\cot(\tilde\theta_{0}(x,\lambda))=\frac{u+Av}{Bv}%
>\cot(D(\lambda)),
\]
or
\[
B(\lambda)\cot(C(\lambda))-A(\lambda)>\frac{u(x,\lambda)}{v(x,\lambda
)}>B(\lambda)\cot(D(\lambda))-A(\lambda),
\]
or
\[
A(\lambda)-B(\lambda)\cot(C(\lambda))<F(x,\lambda)=-\frac{u(x,\lambda
)}{v(x,\lambda)}<A(\lambda)-B(\lambda)\cot(D(\lambda)).
\]
Hence,
\begin{align*}
\mu\left(  \Lambda\cap\left\{  \lambda\in\mathrm{{R\!\!\!\!I} \; \,
|\,C(\lambda)<\tilde\theta_{0}(x,\lambda)<D(\lambda)}\right\}  \right)   &
=\mu(\Lambda\cap\left\{  \lambda\in|\,F(x,\lambda)\in S_{\lambda}\right\}  )
\end{align*}
where $S_{\lambda}:=\left(  A(\lambda)-B(\lambda)\cot(C(\lambda)),A(\lambda
)-B(\lambda)\cot(D(\lambda))\right) $, so it follows from Theorem 2(ii) that
\[
\lim_{x\rightarrow\infty}\mu\left(  \Lambda\cap\left\{  \lambda\in
\mathrm{{R\!\!\!\!I} \; \, |\,C(\lambda)<\tilde\theta_{0}(x,\lambda
)<D(\lambda)}\right\}  \right)  =\int_{\Lambda}\omega(\lambda,S_{\lambda
},m^{+})\,d\lambda.
\]
But, making use of the definitions (\ref{1.3}) and (\ref{1.5}) we observe
that
\begin{align}
\omega(\lambda, S_{\lambda}, m^{+})  & = \frac{1}{\pi} \theta(S_{\lambda},
m^{+})\nonumber\\
& = \frac{1}{\pi} \int_{S_{\lambda}} \operatorname{Im} \left( \frac{1}{t - (A
+ iB)} \right)  \, dt\nonumber\\
& = \frac{1}{\pi} \int_{S_{\lambda}} \frac{-B}{(t-A)^{2} + B^{2}} \,
dt,\nonumber
\end{align}
and making the change of variable $t = A - B \cot\theta$, we obtain
\begin{align}
\omega(\lambda, S_{\lambda}, m^{+})  & = \frac{1}{\pi} \int_{C(\lambda
)}^{D(\lambda)} \frac{B^{2} \csc^{2} \theta}{B^{2} (\cot^{2} \theta+ 1)} \,
dt\nonumber\\
& = \frac{1}{\pi} (D(\lambda) - C(\lambda).\nonumber
\end{align}
Substitution of this into the integral of $\omega(\lambda, S_{\lambda},
m^{+})$ thus yields the result (\ref{2.34}). \quad\rule{2.2mm}{3.2mm}

Theorem 3 gives rise to the concept of a sequence of ``uniformly
asymptotically distributed functions modulo $\pi$''. Namely, when
$g_{n}(\lambda)$ is a sequence of functions defined on $I\subset \sigma_ac$
$\hspace{0.05in}$ and mapping into $[0,\pi]$ (such as $\tilde\theta
_{0}(x,\lambda)$) which satisfy, for any two Lebesgue measurable functions
$C(\lambda),D(\lambda)$ (with $0 \leq C(\lambda) < D(\lambda) \leq\pi$), and
any interval $\Lambda \subset I$, the property (\ref{2.34}), that is,%

\[
\lim_{n\rightarrow\infty}\mu(\Lambda\cap\{\lambda:|\;C(\lambda)<g_{n}%
(\lambda)<D(\lambda)\}=\frac{1}{\pi}\int_{\Lambda}\;\{D(\lambda)-C(\lambda
)\}d\lambda,
\]
we will say that $g_{n}(\lambda)$ is uniformly asymptotically distributed
modulo $\pi$. Further investigations of such functions and how they give rise
to formulas for spectral density functions defined as in (\ref{2.26}) are
currently in progress.

\section{An example:  Bessel's equation}

\setcounter{equation}{0}

{ \fontfamily{times}\selectfont
 \noindent In this section we give an example for which the spectral density function is known,
and show that it can be obtained from the formula in \cite[Thm 4.12]{IMA}.
Consider the Bessel equation of order $\nu \in [0,\infty)$ on $[a,\infty)$, $a>0$,
\bigskip%
\begin{align}
-y"+\left(  \frac{\nu^{2}-1/4}{x^{2}}\right)  y  & =\lambda y,\text{ \ a}\leq
x<\infty \label{3.1}\\
y(a)  & =0. \nonumber
\end{align}
Since this problem has a.c. spectrum in $(0,\infty)$ we may take $I=(0,\infty)$ in Theorem 3.
The solutions of this equation defined by the initial conditions (\ref{2.1}) at $x=a$ are (compare \cite[p. 86]{TITCH}),
\bigskip%
\begin{align}
u(x,\lambda)  & =D_{1}\sqrt{x}J_{\nu}\left(  x\sqrt{\lambda}\right)
+D_{2}\sqrt{x}Y_{\nu}\left(  x\sqrt{\lambda}\right) \label{3.2}  \\
v(x,\lambda)  & =D_{3}\sqrt{x}J_{\nu}\left(  x\sqrt{\lambda}\right)
+D_{4}\sqrt{x}Y_{\nu}\left(  x\sqrt{\lambda}\right), \label{3.3}
\end{align}
where
\bigskip%
\begin{align*}
D_{1}  & =-\frac{\pi}{2}a^{\frac{1}{2}}Y_{\nu}(a\sqrt{\lambda})\text{ \ ,
\ \ D}_{2}=\frac{\pi}{2}a^{\frac{1}{2}}J_{\nu}(a\sqrt{\lambda}),\\
D_{3}  & =\frac{\pi}{4}a^{-\frac{1}{2}}Y_{\nu}\left(  a\sqrt{\lambda}\right)
+\frac{\pi}{2}a^{\frac{1}{2}}\sqrt{\lambda}Y_{\nu}^{\prime}\left(  a\sqrt{\lambda
}\right),  \\
D_{4}  & =-\frac{\pi}{4}a^{-\frac{1}{2}}J_{\nu}\left(  a\sqrt{\lambda}\right)
-\frac{\pi}{2}a^{\frac{1}{2}}\sqrt{\lambda}J_{\nu}^{\prime}\left(
a\sqrt{\lambda}\right).
\end{align*}

\bigskip%
\noindent Using the fundamental system of solutions of the Appell system (\ref{2.18}), we can
find the unique solution, for $\lambda \in (0,\infty)$, of the initial value problem at
infinity (see \cite[Equa (1.6)-(1.7)]{IMA} or \cite[Thm 1]{FPP2}),
\begin{equation} \label{3.4}
\lim_{x\rightarrow\infty} [ \tilde a U_{1}+ \tilde b U_{2}+ \tilde c U_{3} ] =\left(
\begin{array}
[c]{c}%
\sqrt{\lambda}\\
0\\
\frac{1}{\sqrt{\lambda}}%
\end{array}
\right),
\end{equation}
by employing asymptotics for the Bessel functions $J_{\nu}$ and $Y_{\nu}$, as $x\rightarrow\infty$, and this yields
the formulas

\begin{align*}
\tilde a(\lambda)  & =\frac{\pi a}{2}\left(  J_{\nu}^{2}\left(  a\sqrt{\lambda
}\right)  +Y_{\nu}^{2}\left(  a\sqrt{\lambda}\right)  \right),  \\
\tilde b(\lambda)  & =\frac{\pi}{2}\left[  J_{\nu}^{2}\left(  a\sqrt{\lambda}\right)
+Y_{\nu}^{2}\left(  a\sqrt{\lambda}\right)  \right]  +\pi a\sqrt{\lambda
}\left[  J_{\nu}\left(  a\sqrt{\lambda}\right)  J_{\nu}^{\prime}\left(
a\sqrt{\lambda}\right)  +Y_{\nu}(a\sqrt{\lambda})Y_{\nu}^{\prime}\left(
a\sqrt{\lambda}\right)  \right]. \\
\tilde c(\lambda)  & =\frac{\pi}{8a}\left[  J_{\nu}^{2}\left(  a\sqrt{\lambda
}\right)  +Y_{\nu}^{2}\left(  a\sqrt{\lambda}\right)  \right]  +\frac{\pi
a\lambda}{2}\left[  \left(  J_{\nu}^{\prime}(a\sqrt{\lambda}\right)
^{2}+\left(  Y_{\nu}^{\prime}(a\sqrt{\lambda})\right)  ^{2}\right]  \\
& +\frac{\pi\sqrt{\lambda}}{2}\left[  J_{\nu}\left(  a\sqrt{\lambda}\right)
J_{\nu}^{\prime}\left(  a\sqrt{\lambda}\right)  +Y_{\nu}(a\sqrt{\lambda
})Y_{\nu}^{\prime}\left(  a\sqrt{\lambda}\right)  \right].
\end{align*}
From \cite[Thm 4.12, Equa 4.29]{IMA} it follows that (in agreement with \cite[p. 86]{TITCH}),
\begin{equation} \label{3.5}
  f(\lambda) = \frac{1}{\pi \tilde a(\lambda)} = \frac{2}{\pi^2 a} \frac{1}{\left[  J_{\nu}^{2}\left(  a\sqrt{\lambda}\right) +Y_{\nu}^{2}\left(  a\sqrt{\lambda}\right)  \right]}
\end{equation}
From \cite[Equas 4.22 \& 4.26]{IMA} we also have $A(\lambda)+iB(\lambda)= \frac{\tilde b}{2\tilde a} + i\frac{1}{\tilde a}$. Hence
putting these formulas, together with the solutions $\{u,v\}$ from (\ref{3.2})-(\ref{3.3}) into (\ref{2.29}), we have produced,
in the definition (\ref{2.33}) of  $\tilde\theta_{0}(x,\cdot)$, a concrete example of a uniformly asymptotically distributed
sequence of functions modulo $\pi$.
}

 
 \vskip 3mm
 \noindent  {\bf Acknowledgments}


{ \fontfamily{times}\selectfont
 \noindent
The authors thank Jeremy Mandelkern for carrying out the calculations yielding the above formulas for the connection coefficients
${\tilde a,\tilde b, \tilde c}$.}


\begin{thebibliography}{99}                                                                                               

\bibitem {BP2003}S.V. Breimesser and D.B. Pearson,(2003) Geometrical Aspects of
spectral theory and Value Distribution for Herglotz functions. Math Physics,
Anal and Geom, 6, 29--57.

\bibitem {SLEDGE}S. Pruess and C. Fulton, (1993) Mathematical software for
Sturm-Liouville problems. ACM Trans. Math. Software, 19, 360--376.

\bibitem {TOMS98}C. Fulton and S. Pruess, The computation of spectral density
functions for singular Sturm-Liouville problems involving simple continuous
spectra. (1998) ACM Trans. Math. Software, 34, 107--129.

\bibitem {PRYCE}J.D. Pryce, (1993)  Numerical Solution of Sturm-Liouville
Problems, Oxford, Clarendon Press.

\bibitem {FPP2}C. Fulton, D. Pearson, and S. Pruess, (2008) New characterizations of
spectral density functions for singular Sturm-Liouville problems. J. Comp.
Appl. Math, 212, 194--213.

\bibitem {FPP4}C. Fulton, D. Pearson, and S. Pruess, (2008) Efficient calculation of
spectral density functions for specific classes of singular Sturm-Liouville
problems. J. Comp. Appl. Math, 212, 150--178.


\bibitem {IMA}C. Fulton, D. Pearson, and S. Pruess, (2014) Estimating spectral
density functions for Sturm-Liouville problems with two singular endpoints.
IMA Jour. of Numer. Anal. 34 (2014), 609-650.

\bibitem {DBP3}D.B. Pearson, (1993) Value Distribution and spectral analysis of
differential operators. J. Phys. A: Math Gen, 26, 4067--4080.

\bibitem {DBP4}D.B. Pearson,(1994) Value distribution and spectral theory. Proc.
London Math. Soc, 68 (3), 127--144.

\bibitem {BP2000}S.V. Breimesser and D.B. Pearson,(2000) Asymptotic value
distribution for solutions of the Schr\"odinger Equation. Math Physics, Anal
and Geom, 3, 385--403.

\bibitem {AlN2}I. Al-Naggar and D.B. Pearson, (1996) Quadratic forms and solutions of
the Schr\"odinger equation. J. Phys A: Math Gen, 29, 6581-6594.

\bibitem {APPELL}M. Appell, (1880) Sur la transformation des \'{e}quations
diff\'{e}rentielles lin\'{e}aires. Comptes rendus hebdomadaires des
s\'{e}ances de l'Acad\'{e}mie des sciences, 91 (4), 211-214.
                                                                                               
\bibitem {AlN1}I. Al-Naggar and D.B. Pearson, (1994) A new asymptotic condition for
absolutely continuous spectrum of the Sturm-Liouville operator on the
half-line. Helvetica Physica Acta 67, 144-166.

\bibitem{TITCH}Titchmarsh, E.C., (1962) Eigenfunction Expansions associated with
Second-order Differential Equations, 2nd edn., Oxford, Clarendon Press. 




\end{thebibliography}
\end{document}